\documentclass[12pt, twoside, leqno]{article}

\usepackage{amsfonts}
\usepackage{amsthm}
\usepackage{amsmath}
\usepackage{amsxtra}
\usepackage{amstext}
\usepackage{amssymb}
\usepackage{latexsym}
\numberwithin{equation}{section}
\usepackage{epsfig}
\usepackage{graphicx}
\usepackage[subnum]{cases}
\usepackage[usenames,dvipsnames]{color}
\usepackage{enumerate}

\usepackage[english]{babel}
\usepackage[ansinew]{inputenc}

\theoremstyle{definition}
\newtheorem{fed}{Definition}[section]

\newtheorem{teo}[fed]{Theorem}
\newtheorem*{teo*}{Theorem}

\newtheorem{lem}[fed]{Lemma}
\newtheorem{cor}[fed]{Corollary}
\newtheorem{pro}[fed]{Proposition}
\newtheorem{obs}[fed]{Remark}

\newtheorem{ejem}[fed]{Example}

\newtheorem{ej}{Ejercicio}[subsection]

\numberwithin{equation}{section}

\def\bteo{\begin{teo}}
\def\eteo{\end{teo}}
\def\blem{\begin{lem}}
\def\elem{\end{lem}}
\def\bpro{\begin{pro}}
\def\epro{\end{pro}}
\def\bcor{\begin{cor}}
\def\ecor{\end{cor}}
\def\bobs{\begin{obs}}
\def\eobs{\end{obs}}
\def\bdefi{\begin{fed}}
\def\edefi{\end{fed}}
\def\bdem{\begin{proof}}
\def\edem{\end{proof}}
\def\bexa{\begin{ejem}}
\def\eexa{\end{ejem}}
\def\bexe{\begin{ej}}
\def\eexe{\end{ej}}

\def\r{\mathbb{R}}
\def\R{\mathbb{R}}

\def\l{\mathcal L}
\def\M{\mathcal M}
\def\m{\mathcal M}

\def\O{\Omega}

\begin{document}
\baselineskip=17pt

\title{Weighted a priori estimates for elliptic equations}

\author{Mar\'\i a E. Cejas \\
	Departamento de Matem\'atica \\
	Facultad de Ciencias Exactas \\
	Universidad Nacional de La Plata \\ 
	CONICET \\
	Calle 50 y 115,	1900 La Plata, Buenos Aires, Argentina \\
	E-mail: mec.eugenia@gmail.com
\and
	Ricardo G. Dur\'an \\
	Departamento de Matem\'atica\\ 
	Facultad de Ciencias Exactas y Naturales \\
	Universidad de Buenos Aires\\
	IMAS-CONICET-UBA \\
	Pabell\'on I, Ciudad Universitaria, 1428, CABA, Argentina \\
	E-mail:rduran@dm.uba.ar}
\date{}
\maketitle

\renewcommand{\thefootnote}{}



\renewcommand{\thefootnote}{\arabic{footnote}}
\setcounter{footnote}{0}

\begin{abstract}
We give a simpler proof of the a priori estimates obtained in
\cite{DST1,DST2} for solutions of elliptic problems in weighted Sobolev
norms when the weights belong to the Muckenhoupt class $A_p$.
The argument is a generalization to bounded domains of the one used
in $\R^n$ to prove the continuity of singular integral operators
in weighted norms.

In the case of singular integral operators it is known that
the $A_p$ condition is also necessary for the continuity. We do not
know whether this is also true for the a priori estimates in bounded domains
but we are able to prove a weaker result when the operator is the Laplacian
or a power of it. We prove that a necessary condition is
that the weight belongs to the local $A_p$ class.
\end{abstract}

\footnotetext[3]{}

\section{Introduction}
The goal of this paper is to prove weighted a priori estimates for
solutions of linear elliptic problems with Dirichlet boundary conditions.
More precisely, for a bounded smooth domain $\Omega\subset\R^n$ we consider

$$
\left\{
\begin{aligned}
\l u&=f   \mbox{ in } \Omega \\
\mathcal{B}_ju&=0 \mbox{ on } \partial \Omega \hspace{0,2cm},\quad 1 \leq j \le m-1
\end{aligned}
\right.
$$
where $\l$ is an elliptic operator of order $2m$ and $\mathcal{B}_j$ differential
operators of order $m_j$ satisfying the properties introduced in the
classic paper \cite{ADN}.

For $1<p<\infty$, the a priori estimate
\begin{equation}
\label{apriori adn}
\|u\|_{W^{2m,p}(\Omega)} \le C \|f\|_{L^p(\Omega)}.
\end{equation}
is well-known.
This result is usually referred as Agmon-Douglis-Nirenberg
estimate because it is essentially contained in \cite{ADN} although
it is not explicitly written in this way in that paper. For completeness
we give more details on this point in the next section.

The estimate \eqref{apriori adn} has been extended to weighted norms
when $\l$ is the Laplacian or a power of it
in \cite{DST1,DST2}.

By a weight function we mean a locally integrable function $w$ defined in
$\R^n$ and for $1\le p<\infty$, we define the Banach space $L_w^p(\Omega)$
with norm given by
$$
\|f\|_{L_w^p(\Omega)}=\left(\int_{\Omega}|f(x)|^pw(x)\,dx\right)^{1/p}.
$$
and $W_w^{2m,p}(\Omega)$ the Sobolev
space of functions in $L^p_w(\Omega)$ with derivatives up to order $2m$
in $L^p_w(\Omega)$ with the usual norm. For $\O=\R^n$ we write simply $L^p_w$ instead of $L^p_w(\R^n)$.
With $C$ we will denote a generic constant which can change its value
even in the same line.

For $1<p<\infty$, a weight is in the Muckenhoupt class $A_p$ if
\begin{equation}
\label{muckenhoupt}
\left(\frac{1}{|Q|}\int_Q w\right)
\left(\frac{1}{|Q|}\int_Q w^{-\frac{1}{p-1}}\right)^{p-1} \le C
\end{equation}
for all cubes $Q$. It is well known that the Hardy-Littlewood maximal operator is
bounded in $L^p_w$ if and only if $w\in A_p$ (see for example \cite{D}).

For $\l=(-\Delta)^m$ it was proved in \cite{DST2} (extending the results for $m=1$
given in \cite{DST1}) that, if $w$ is in $A_p$, then
\begin{equation}
\label{estimacion a priori}
\|u\|_{W_w^{2m,p}(\Omega)} \le C \|f\|_{L_w^p(\Omega)},
\end{equation}

In this paper we give a different proof of \eqref{estimacion a priori}
generalizing to a bounded domain
the classic arguments used to obtain the continuity of singular integral
operators in $L_w^p$. The advantage of this proof is that it is simpler
and it does not require estimates of derivatives of the Green function
involving the distance to the boundary. On the other hand, our arguments are
very general and apply to the class of operators considered in \cite{ADN,K}.

We do not know whether the $A_p$ condition is also necessary to have
(\ref{estimacion a priori}) but we prove a weaker result
for the case of $\l=(-\Delta)^m$. Indeed, in order to
have the a priori estimate (\ref{estimacion a priori}) it is necessary
that $w\in A^{loc}_p(\O)$ (see below for the definition of this class).

\section{Weighted a priori estimates}

For a bounded domain $\O\subset\R^n$ we consider the problem
\begin{equation}
\label{eliptico}
\begin{aligned}
\l u(x):=&\sum_{0\le |\alpha|\le 2m} a_{\alpha}(x)D^{\alpha}u(x)=f(x)
\quad\mbox{ in } \Omega \\
\mathcal{B}_ju(x):=&\sum_{0\le|\alpha|\le m_j}  b_{j\alpha}(x) D^{\alpha}u(x)=0
\quad\quad\mbox{ on } \partial\Omega \hspace{0,2cm} 1 \le j \le m-1
\end{aligned}
\end{equation}
where $\l$ is uniformly elliptic and ${\mathcal B}_j$ satisfy the complementing
conditions in the sense of \cite{ADN}.

Our arguments are based on the estimates for the Green function of \eqref{eliptico}
proved in \cite{K}, and therefore, we assume the hypotheses of that paper.
Namely, $m_j\le 2m-1$ and,
for $\ell_0:=\max_j (2m-m_j)$, let $\ell_1$ be an integer such that
\begin{equation}
\label{regularidad}
\begin{aligned}
&\ell_1\ge \frac32 \ell_0 , \quad \mbox{for} \ \ n\ge 3\\
&\ell_1\ge 2(\ell_0+1) , \quad \mbox{for} \ \ n=2.
\end{aligned}
\end{equation}
Then,
\begin{equation}
\label{coeficientes}
a_\alpha\in C^{\ell_1+1}(\O) \qquad  \qquad
b_{j\alpha}\in C^{\ell_1+1}(\partial\O)
\end{equation}
and
\begin{equation}
\label{regularidad del borde}
\partial\O\in C^{\ell_1+2m+1}
\end{equation}
Under these assumptions, the existence of the Green function as well as
estimates for it and its derivatives were proved in \cite{K}. We state in the next
theorem the estimates that we are going to use in our arguments.

\bteo
\label{estimaciones green}
Under the hypotheses (\ref{regularidad}), (\ref{coeficientes})
and (\ref{regularidad del borde}) there exists the Green function $G$
of \eqref{eliptico}, namely, the solution $u$ is given by
\begin{equation}
\label{solucion}
u(x)=\int_{\Omega} G(x,y)f(y)\,dy .
\end{equation}
Moreover, for $0<\alpha<1$,
there exists a constant $C$ depending on $\l$, ${\mathcal B}_j$, $\Omega$ and $n$ such that

\begin{equation}
\label{cotas Green 1}
|\gamma|\le 2m \quad , \quad 2m-n-|\gamma|\neq 0
\quad \Longrightarrow \quad
|D_x^{\gamma}G(x,y)|\le C \left\{|x-y|^{2m-n-|\gamma|}+1\right\}
\end{equation}

\begin{equation}
\label{cotas Green 2}
\quad 2m-n-|\gamma|=0
\quad \Longrightarrow \quad
|D_x^{\gamma}G(x,y)| \le C \left\{|\log|x-y||+1\right\}
\end{equation}

\begin{equation}
\label{cotas Green 3}
|\gamma|=2m
\quad \Longrightarrow \quad
\left|D_x^{\gamma}G(y,z)-D_x^{\gamma}G(x,z)\right|
\leq C \left|y-x\right|^{\alpha}\left(\left|y-z\right|^{-n-\alpha}+\left|x-z\right|^{-n-\alpha}\right).
\end{equation}
\eteo
\bdem See Theorem 3.3 and its Corollary in \cite[Page 965]{K}. \edem

Our proof of the weighted a priori estimates makes use of the classic
unweighted results. As we mentioned in the introduction, the a priori estimate
\eqref{apriori adn} is essentially contained in \cite{ADN} although not
written explicitly there. Indeed, in the particular case of homogeneous boundary conditions
that we are considering, Theorem 15.2 in page 704 of \cite{ADN} says that
$$
\label{apriori adn verdadera}
\|u\|_{W^{2m,p}(\Omega)} \le C \left\{\|f\|_{L^p(\Omega)} + \|u\|_{L^p(\Omega)}\right\}.
$$

But, in view of the representation \eqref{solucion} and the bound for $G$ given in
\eqref{cotas Green 1}, a standard application of the Young inequality
yields $\|u\|_{L^p(\Omega)}\le C\|f\|_{L^p(\Omega)}$, and therefore, \eqref{apriori adn}
follows.

Let us remark that in the above mentioned Theorem 15.2 of \cite{ADN} the authors
assume that the norm on the left hand side of the estimate is finite, but this follows from
their previous Theorem 7.3 (see also Remark 1 after that theorem) \cite[Page 668]{ADN}
provided $f$ is regular enough and using again the representation \eqref{solucion}
to bound the norm of $u$ appearing in the right hand side of that theorem.
Then, in the general case one can proceed by a standard density argument.

Let us now recall the argument used in the case of
singular integral operators that we are going to generalize.
We will make use of the Hardy-Littlewood maximal operator
$$
\M f(x)= \sup_{Q \ni x} \frac{1}{|Q|} \int_Q  |f(y)|\,dy
$$
and of the sharp maximal operator
$$
\M^{\#}f(x)= \sup_{Q \ni x }\frac{1}{|Q|}
\int_Q |f(y)-f_Q|\,dy,
$$
where the supremums are taken over all cubes containing $x$ and $f_Q:=\frac1{|Q|}\int_Qf$.

If
$$
Tf(x)=\lim_{\varepsilon\to 0}\int_{|x-y|>\varepsilon}K(x,y)f(y)dy
$$
is a singular integral operator which is continuous in $L^p$, for
$1<p<\infty$, and $K(x,y)$ satisfies
$$
|K(x,z)-K(y,z)|\le \frac{C|x-y|^\alpha}{|x-z|^{n+\alpha}}\ ,
\qquad \mbox{for} \quad |x-z|\ge 2
|x-y|
$$
with $0<\alpha<1$, then we have, for any $s>1$,
\begin{equation}
\label{cota M sharp}
\M^\#Tf(x)\le C (\M|f|^s(x))^{1/s}.
\end{equation}
This
estimate is well-known and its proof can be found in several
books, although the hypotheses on the operator are not stated
usually as we are doing here. Indeed, the proof given in
\cite[Lemma 7.9]{D} only uses the hypotheses given above.

In the next lemma we prove a version of \eqref{cota M sharp} in a bounded domain.
With this goal we introduce the local
sharp maximal operator
$$
\M_{\Omega}^{\#}f(x)= \sup_{ \Omega \supset Q \ni x }\frac{1}{|Q|}
\int_Q |f(y)-f_Q|\,dy.
$$

\blem
\label{puntual}
Let $u$ be the solution of \eqref{eliptico} and assume that the hypotheses
of Theorem \ref{estimaciones green} are satisfied.
If $|\gamma|=2m$ we have, for any $s>1$ and any $x\in\O$,
$$
\M^\#_{\Omega}(D^\gamma u)(x)
\le  C (\m\left|f\right|^s)^{1/s}(x).
$$
\elem
\bdem We extend $f$ by zero outside of $\O$.
Let $Q \subset \Omega$ be a cube such that $x\in Q$ and $Q^*$ an expanded cube of $Q$
by a factor $2$. We decompose $f=f_1+f_2$ with $f_1=f \chi_{Q^*}$ and call
$u_1$ and $u_2$ the solutions of \eqref{eliptico} with $f_1$ and $f_2$ as right hand sides
respectively.

It is known that we can replace $f_Q$ by any constant.
We choose $a=D^\gamma u_2(x)$.  Then,
$$
\frac{1}{\left|Q\right|} \int_{Q} \left|D^\gamma u(y)-a\right|\,dy
\le  \frac{1}{\left|Q\right|} \int_{Q} \left|D^\gamma u_1(y)\right|\,dy
+\frac{1}{\left|Q\right|} \int_{Q} \left|D^\gamma u_2(y)-D^\gamma u_2(x)\right|\,dy
= (i)+(ii).
$$
Now, given $s>1$, using the H\"older inequality and \eqref{apriori adn},
we have

$$
\begin{aligned}
(i)\le\left(\frac{1}{\left|Q\right|}
\int_Q\left|D^\gamma u_1(y)\right|^s \,dy\right)^{1/s}
&\le C\left(\frac{1}{\left|Q\right|}\int_\O
\left|f_1(y)\right|^s \,dy \right)^{1/s} \\
= C \left(\frac{1}{|Q|}
\int_{Q^*} \left|f(y)\right|^s\,dy\right)^{1/s}
&\le C (\m \left|f\right|^s)^{1/s}(x).
\end{aligned}
$$
On the other hand, if $x\notin{\rm supp\,}f_2$ we can take the derivative
inside the integral in the expression for $u_2$ given by \eqref{solucion}.
Then, since ${\rm supp\,}f_2\subset (Q^*)^c$, for $x\in Q$
we have
$$
(ii)
\le\frac{1}{\left|Q\right|}\int_{Q}\int_{(Q^*)^c}
\left|D_x^{\gamma}G(y,z)-D_x^{\gamma}G(x,z)\right|\left|f_2(z)\right|\,dz\,dy.
$$
Therefore, using \eqref{cotas Green 3} and that
$\left|y-z\right|\sim \left|x-z\right|$
and $\left|x-z\right|\geq \frac{\ell(Q)}{2}$ we obtain
$$
\begin{aligned}
(ii)
&\le C\frac{l(Q)^\alpha}{\left|Q\right|}\int_{Q}\int_{(Q^*)^c}
\frac{|f(z)|}{|x-z|^{n+\alpha}}\,dz\,dy
\le C l(Q)^{\alpha}\int_{(Q^*)^c} \frac{|f(z)|}{|x-z|^{n+\alpha}} \,dz \\
&\le
C l(Q)^{\alpha}\sum_{k=0}^\infty
\int_{2^{k-1}l(Q)<|z-x|\le 2^kl(Q)}\frac{|f(z)|}{|x-z|^{n+\alpha}} \,dz
\le C \M f(x)
\end{aligned}
$$
where the last inequality follows by a standard argument
(see \cite[Lemma 7.9]{D} for details).
\edem

The following Lemma is a slightly modified version of \cite[Theorem 5.23]{DRS} because
we are using a different definition of the sharp maximal operator. The reader can easily check
that the proof given in that paper applies to our case.

\blem
\label{Fefferman-Stein}
For $f\in L^1_{loc}(\O)$
and $w\in A_p$, if  $f_{\Omega}=\frac{1}{|\Omega|} \int_{\Omega}|f|$, then
\begin{equation}
\label{desigualdad fefferman}
\|f-f_{\Omega}\|_{L_w^p(\Omega)}
\le C \|\m_{ \Omega}^{\#}f\|_{L_w^p(\Omega)}.
\end{equation}
\elem
Now we can state and prove the main result of this section.

\bteo
\label{teo a priori}
Let $u$ be the solution of \eqref{eliptico} and assume that the hypotheses
of Theorem \ref{estimaciones green} are satisfied. Then, for $1<p<\infty$,
$w\in A_p$ and $f\in L_w^p(\Omega)$, there
exists a constant $C$ depending on $\l$, ${\mathcal B}_j$, $\O$, $n$ and $w$ such that
$$
\|u\|_{W_w^{2m,p}(\Omega)}
\leq C \left\|f\right\|_{L_w^p(\Omega)}.
$$
\eteo
\bdem If $2m-n<|\gamma|<2m$, using \eqref{solucion} and
\eqref{cotas Green 1} we obtain
$$
\begin{aligned}
|D^\gamma u(x)|&=\left|\int_\O D_x^\gamma G(x,y) f(y)\,dy\right|
\le C \int_\O\frac{|f(y)|}{|x-y|^{|\gamma|+n-2m}}\,dy
\\
&\le C \sum_{k=0}^\infty
\int_{2^{-(k+1)}d<|x-y|\le 2^{-k}d} \frac{|f(y)|}{|x-y|^{|\gamma|+n-2m}}\,dy
\le C\M f(x)
\end{aligned}
$$
where $d$ denotes the diameter of $\O$.
For $|\gamma|\le 2m-n$ we obtain the same estimate using now that,
in view\eqref{cotas Green 1} and
\eqref{cotas Green 2}, for any $\varepsilon>0$,
$|G(x,y)|\le C |x-y|^{-\varepsilon}$.

Consequently, it follows from the boundedness of the maximal operator for
$A_p$ weights that
$$
\|u\|_{W_w^{2m-1,p}(\Omega)}
\le C \left\|f\right\|_{L_w^p(\Omega)}.
$$
It rests to estimate $D^\gamma u$ for
$|\gamma|=2m$ which is the most difficult part.
From Lemmas \ref{Fefferman-Stein} and \ref{puntual} we have

\begin{equation}
\label{cota T gama 1}
\begin{aligned}
\int_{\Omega} \left|D^\gamma u(x)- (D^\gamma u)_\O\right|^p w(x)\,dx
&\le C \int_\O \left|\M_\O^{\#}(D^\gamma u)(x)\right|^p w(x)\,dx \\
&\le C\int_\O (\M\left|f\right|^s(x))^{p/s}w(x)\,dx.
\end{aligned}
\end{equation}
But it is known that there exists $s$
depending only on $w$ such that $1<s<p$ and $w \in A_{p/s}$ (see for example \cite{D}),
and using the boundedness of $\M$ in $L^{p/s}_w$ we obtain from \eqref{cota T gama 1},
$$
\int_{\Omega} \left|D^\gamma u(x)- (D^\gamma u)_\O\right|^p w(x)\,dx
\le C \int_\O \left|f(x)\right|^p w(x)\,dx.
$$
Then,
$$
\int_{\Omega} \left|D^\gamma u(x)\right|^p w(x)\,dx
\le C \left(\int_\O\left|f(x)\right|^p w(x)\,dx + \int_\O |(D^\gamma u)_\O|^p w(x) \,dx\right)
$$
and so it only remains to estimate the last term. But,
since $w$ is integrable in $\Omega$, it is enough to show that
\begin{equation}
\label{cota del promedio}
\left|(D^\gamma u)_\O\right|\le C \left\|f\right\|_{L_w^p(\Omega)}.
\end{equation}
Taking $1<s<p$ and using the a priori estimate \eqref{apriori adn} for $s$
we have
$$
\begin{aligned}
\left|(D^\gamma u)_\O\right|
\le
\left(\frac{1}{|\O|}\int_{\Omega} \left|(D^\gamma u)_\O\right|^s \,dx \right)^\frac{1}{s}
&\le C \left(\frac{1}{|\O|}\int_{\Omega} |f(x)|^s \,dx \right)^\frac{1}{s}
\\
&\le  C \left(\frac{1}{|\O|}\int_{\Omega} |f(x)|^p w(x) \,dx\right)^\frac{1}{p}
\left(\frac{1}{|\O|}\int_{\Omega} w(x)^{-\frac{s}{p-s}}\,dx \right)^{\frac{p-s}{sp}}
\end{aligned}
$$
where we have used the H\"older inequality first with $s$
and then with $p/s$. Then, choosing $s$ such that $w \in A_{p/s}$ the last term on the right
hand side is finite,
and therefore, \eqref{cota del promedio} is proved.
\edem

\section{Necessary condition}

It is known that the $A_p$ condition is also necessary for the continuity of
singular integral operators \cite{S}. Then, it is natural to ask whether the
same is true for the weighted a priori estimates. We do not know the answer
but we prove in this section a weaker result, namely, a necessary condition to
have the weighted a priori estimates for $\l=(-\Delta)^m$ is that the weight belong
to the $A_p^{loc}(\O)$ class extensively studied in \cite{HSV}.

In this section it is more convenient to work with balls instead of cubes.
To recall the definition of the $A_p^{loc}$ class, first we consider $0<\beta<1$ and define the
$A_p^{\beta}(\O)$ class as follows. A weight $w$ belongs to this class if
$$
\left(\frac{1}{|B|}\int_B w \right)\left(\frac{1}{|B|}\int_B w^{-\frac{1}{p-1}} \right)^{p-1}
\le C
$$
for all balls $B\subset\O$ such that $diam(B)<\beta dist(B,\partial \Omega)$. It was proved in
\cite{HSV} that the classes $A _p^{\beta}(\Omega)$ are independent of $\beta$,
namely, if $0<\beta<1$ and $0<\gamma<1$, we have that $w \in A_p^{\beta}(\Omega)$
if and only if
$w \in  A_p^{\gamma}(\Omega)$.
In view of this fact, we say that $w \in A_p^{loc}(\Omega)$ if $w \in A_p^{\beta}(\Omega)$
for some $0<\beta<1$.
We will call a ball satisfying this condition for some $\beta$, admissible.
To simplify notation, we will use the usual notation
$w(S)=\int_S w(x)\,dx$.

\bpro
\label{condicion equivalente}
Let $w$ be a weight. Then, $w \in A_p^{loc}(\Omega)$ if and only if
\begin{equation}
\label{equivalencia aploc}
(f_B)^p \leq \frac{C}{w(B)}\int_B f^p w \,dx
\end{equation}
for all $f$ nonnegative and for all admissible ball $B$.
\epro
\bdem The proof follows as in the case of the $A_p$ weights given
in \cite[Chapter V, Section 1.4]{S} and \cite[Theorem 7.1]{D} using the results in \cite{HSV}
and considering admissible balls.
\edem
Let $f \in L_w^p(\Omega)$ and consider the homogeneous boundary value problem

\begin{equation}
\label{potencia laplaciano}
\left\{
\begin{aligned}
(-\Delta)^m u&=f \quad \mbox{ in } \Omega \\
\left(\frac{\partial}{\partial \nu}\right)^{j}u&=0 \quad\mbox{ on }
\partial \Omega \hspace{0,2cm},\quad 1 \le j \le m-1
\end{aligned}
\right.
\end{equation}
then
\begin{equation}
\label{representacion lap a la m}
u(x)=\int_{\Omega}G_m(x,y)f(y)\,dy
\end{equation}
where $G_m(x,y)$ is the Green function of the operator $(-\Delta)^m$ in $\Omega$. It is well
known that
$$
G_m(x,y)= \Gamma(x-y)+h(x,y)
$$
where $\Gamma$ is the fundamental solution given by
$$
\left\{
\begin{aligned}
&c_{m,n}|x|^{2m-n} \qquad & n\ \mbox{odd, or}
\ n\ \mbox{even and}\ n > 2m
\\
&c_{m,n}|x|^{2m-n}\log |x| \qquad & n  \mbox{ even and } n \le 2m
\end{aligned}
\right.
$$
and, for each $y\in \Omega$, $h(x,y)$ satisfies
$$
\left\{
\begin{aligned}
(-\Delta_x)^m h(x,y)&=0 \hspace{3.7cm} x \in \Omega \\
\\
\left(\frac{\partial}{\partial \nu}\right)^{j}h(x,y)
&=- \left(\frac{\partial}{\partial \nu}\right)^{j}\Gamma(x-y)
\qquad x \in \partial\Omega\ ,\ 0 \le j \le m-1.
\end{aligned}
\right.
$$

Let us recall that, for $|\gamma|=2m$, a standard argument yields
\begin{equation}
\label{parte de gamma}
D_x^{\gamma}\int_\O\Gamma(x-y) f(y)\,dy
=\lim_{\varepsilon \rightarrow 0}
\int_{|x-y|>\epsilon} D^{\gamma}\Gamma(x-y)f(y)\,dy+ c(x)f(x)
\end{equation}
where $c$ is a bounded function.

We will use the ideas given in \cite[Chapter V, Section 4.6]{S}.

\blem
\label{prop fundamental}
For $|\gamma|=2m$ we have
\begin{enumerate}
\item There exists $u_0 \in \r^n$ with $|u_0|=1$ and a constant $C_0$
such that, for all positive numbers $t$,
$$
|D^{\gamma}\Gamma(tu_0)| \ge C_0 t^{-n}.
$$
\item There exists $t_0$ such that if $u=t_0u_0$ and $|v|\leq 2$ then for all $0\neq r\in\R$,
$$
|D^{\gamma}\Gamma(r(u+v))-D^{\gamma}\Gamma(ru)|\leq \frac{1}{2}|D^{\gamma}\Gamma(ru)|.
$$
\end{enumerate}
\elem

\bdem
One can check that $D^{\gamma}\Gamma$ is homogeneous of degree $-n$
and not identically zero. Then,
there exists $u_0 \in \r^n$ with $|u_0|=1$ such that $|D^{\gamma}\Gamma(u_0)|=:C_0>0$,
and therefore,
(1) follows from the homogeneity.

To prove (2) we observe first that, by homogeneity, it is enough to show that
the statement holds for $r=1$.

Take $v\in\R^n$ satisfying $|v|\le k|u|$
with $k\le \frac12$ to be chosen below. Then, for some $\xi$ in the segment
joining $u$ and $u+v$,
$$
|D^{\gamma} \Gamma(u+v)-D^{\gamma}\Gamma(u)|
\le|\nabla D^{\gamma} \Gamma(\xi)||v|
\le C |\xi|^{-n-1}|v|
$$
and using $|u|\le 2|\xi|$, $|v|\le k|u|$ and (1) we obtain,
$$
|D^{\gamma} \Gamma(u+v)-D^{\gamma}\Gamma(u)|
\le  C_1 k|u|^{-n}
\le \frac{C_1 k}{C_0} |D^{\gamma}\Gamma(u)|.
$$
Consequently it is enough to choose
$k$ such that $\frac{C_1}{C_0} k \le\frac12$.
Now, since $|v|\le 2$, if we choose $t_0=\frac2{k}$ our hypothesis
$|v|\le k |u|$ is verified and the proof is concluded.
\edem

The following result is proved in (\cite[Proposition 3.3]{DST2}).

\blem
\label{prop h}
There exists a constant $C$ such that,
for $|x-y| \le d(x)$,
$$
|D_x^{\gamma}h(x,y)|\le Cd(x)^{-n}.
$$
\elem
Now, using the representation \eqref{representacion lap a la m} and \eqref{parte de gamma},
we have
\begin{equation}
\label{derivada = T gamma}
D^\gamma u(x)= T_\gamma f(x) + c(x) u(x)
\end{equation}
where $T_\gamma$ is defined by
$$
T_\gamma f(x)= \lim_{\varepsilon \rightarrow 0}
\int_{|x-y|>\epsilon} D^{\gamma}\Gamma(x-y)f(y)\,dy
+ \int_\O D_x^{\gamma} h(x,y)f(y)\,dy .
$$
We can now prove the main result of this section.

\bteo
\label{cond necesaria}
Let $u$ be the solution of \eqref{potencia laplaciano}. If
$w$ is a weight such that the following a priori estimate
$$
\|u\|_{W^{2m,p}_w(\Omega)} \leq C \|f\|_{L_w^p(\Omega)}
$$
holds, then $w \in A_p^{loc}(\Omega)$.
\eteo
\bdem
In view of \eqref{derivada = T gamma}, since $c$ is a bounded function, it is enough to prove that,
if for any $\gamma$ such that $|\gamma|=2m$,
$$
\|T_\gamma f\|_{L_w^p(\Omega)}\le C\|f\|_{L_w^p(\Omega)}
$$
then $w \in A_p^{loc}(\Omega)$.

We write $T_\gamma=T_1+T_2$ where
$$
T_1f(x)=\lim_{\varepsilon \rightarrow 0}
\int_{|x-y|>\epsilon} D^{\gamma}\Gamma(x-y)f(y)\,dy
$$
and
$$
T_2f(x)=\int_\O D_x^{\gamma} h(x,y)f(y)\,dy.
$$
Consider an admissible ball $B=B(\bar{x},r)$ , i.e., $2r<\beta dist(B,\partial \Omega)$,
with $\beta$ to be determined later and $B':=B(\bar{x}+ru,r)$, with  $u=t_0 u_0$ as
in Lemma \ref{prop fundamental}. We will see that $\beta$ can be taken so that $B'$
is also admissible. Let $x \in B'$ and $z \in \partial \Omega$ then $x=\bar{x}+ru+rx'$
with $|x'|\le 1$,
$$
\begin{aligned}
|x-z|&\ge |\bar{x}-z|-r|u+x'|
\\&\ge dist(B,\partial \Omega)-r|u+x'|
\\&\ge\frac{2r}{\beta}-r|u+x'|
\\&\ge 2r\left(\frac{1}{\beta}-\frac{t_0+1}{2}\right),
\end{aligned}
$$
so, taking $\beta$ satisfying $\frac{1}{\beta}-\frac{(t_0+1)}{2}>1$, i.e.
$\beta<\frac{2}{t_0+3}$ and $\frac{1}{\alpha}=\frac{1}{\beta}-\frac{(t_0+1)}{2}$,
we have $2r<\alpha dist(B',\partial \Omega)$.

Now, we will show that for $x \in B'$ and $y\in B$, $D^{\gamma}\Gamma(x-y)$ has constant sign.
Indeed, writing $x=\bar{x}+ru+rx'$ and $y=\bar{x}+ry'$ with $|x'|,|y'|\le 1$,
we have $x-y=ru+rv$ with $|v|=|x'-y'|\le 2$. Then, by (2) from Lemma \ref{prop fundamental}
we obtain, for $D^{\gamma}\Gamma(ru)>0$,
\begin{equation}
\label{desi 1}
\frac12 D^{\gamma}\Gamma(ru)\le D^{\gamma}\Gamma(x-y)\leq \frac{3}{2} D^{\gamma}\Gamma(ru),
\end{equation}
while for $D^{\gamma}\Gamma(ru)<0$
\begin{equation}
\label{desi 2}
\frac32 D^{\gamma}\Gamma(ru)\le D^{\gamma}\Gamma(x-y)\leq \frac{1}{2} D^{\gamma}\Gamma(ru).
\end{equation}
Consequently, taking $f \in C_0^{\infty}(B)$ positive, we have
$$
|T_1f(x)|=\int_B |D^{\gamma}\Gamma(x-y)|f(y) \,dy.
$$
and moreover, using (\ref{desi 1}), (\ref{desi 2}) and property $(1)$
of Lemma \ref{prop fundamental},
$$
|T_1f(x)| \ge \frac12 \int_B |D^{\gamma}\Gamma(ru)|f(y)\,dy
\ge C_0(rt_0)^{-n} \int_B f(y)\,dy
= C_1f_B.
$$
with a constant $C_1$ depending only on $t_0$, $C_0$ and $n$.

On the other hand, in order to bound $|T_2f(x)|$ we use Lemma \ref{prop h}.
We require $|x-y|\le d(x)$, but,
$$
|x-y|=|ru+r(x'-y')|\leq r(|u|+|x'-y'|) \leq r(t_0+2) <\frac{\alpha}{2} d(x) (t_0+2)
$$
and then, we need $\frac{\alpha}{2} (t_0+2)< 1$ or equivalently
$\beta<\frac{2}{2t_0+3}$. Now, we have that
$$
\begin{aligned}
|T_2f(x)|\le \int_{B}|D_x^{\gamma}h(x,y)|f(y)\,dy
\le C\int_B d(x)^{-n}f(y)\,dy&\le C (2r)^{-n}\alpha^n \int_Bf(y)\,dy
\\
&= C_2 \alpha^n f_B.
\end{aligned}
$$
where $C$ is the constant appearing in Lemma \ref{prop h} and $C_2$ depends on $n$ and $C$.

Summing up we have
$$|T_\gamma f(x)| \geq (C_1-C_2\alpha^n)f_B.$$

In order to have $C_1-C_2 \alpha^n>0$ it is enough that
$\beta< \frac{1}{\frac{t_0+1}{2}+(\frac{C_2}{C_1})^{\frac{1}{n}}}$.

Taking into account the other conditions for $\beta$ we choose
$\beta<\min \left\{\frac{2}{2t_0+3},
\frac{1}{\frac{t_0+1}{2}+(\frac{C_2}{C_1})^{\frac{1}{n}}}\right\}$,
and then, we obtain $f_B\le C|T_\gamma f(x)|$,
for any $x\in B'$. Therefore,
$$
\int_{B'} f_B^p  w(x)\,dx \le C\int_{B'} |T_\gamma f(x)|^pw(x)\, dx
\le C\int_{\Omega} |T_\gamma f(x)|^p w(x) \,dx
$$
and consequently, applying the continuity of $T_\gamma$ in $L_w^p$,
we obtain
\begin{equation}
\label{eq 1}
(f_B)^pw(B') \le C \int_{B} f(x)^p w(x)\,dx.
\end{equation}
Analogously, changing roles of $B$ and $B'$ and taking $f$ with support in $B'$, it follows
that
\begin{equation}
\label{eq 2} (f_{B'})^p w(B) \leq C \int_{B'} f(x)^p w(x)\,dx.
\end{equation}
provided that $C_1-C_2 \beta^n>0$. But this inequality holds for our
previous election of $\beta$ because $\beta<\alpha$.

By a passage to the limit, the inequalities \eqref{eq 1} and \eqref{eq 2}
extend to any non-negative function $f$  supported in $B$ or $B'$ respectively.
If we consider $f=\chi_{B'}$ in (\ref{eq 2}) we obtain
$$
w(B) \leq C w(B').
$$
Using this in (\ref{eq 1}) we get
$$(f_B)^p w(B) \leq C \int_{B} f(x)^p w(x)\,dx$$
but, by Proposition \ref{condicion equivalente}, this means that
$w \in A_p^{loc}(\Omega)$.\edem

\subsection*{Acknowledgements}
This research was supported by Universidad de
Buenos Aires under grant 20020120100050 and by ANPCyT, Argentina,
under grant PICT 2014-1771.
The results of this paper are part of the
doctoral thesis of the first author which was done under a fellowship granted by
CONICET, Argentina. We thank the anonymous referee for helpful comments
that led us to improve the previous version of this work.

\end{document}